\input amstex
\documentstyle{amsppt}
\vsize8.9 true in
\hsize6.5 true in
\loadeusm
\def\script#1{{\fam\eusmfam\relax#1}}
\input amssym
\nologo

\define\pd{positive definite  }
\define\s{\frak S (\tau)}

\define\LP{Laguerre-P\'olya }
\def\lp{\script{L}\text{-}\script{P}}

\def\Re{\mathop{\text{\rm Re}}}

\def\op{\operatorname}
\def\lra{\longrightarrow}
\def\ds{\displaystyle}
\def\re{\mathop{\text{\rm Re}}}
\def\im{\mathop{\text{\rm Im}}}
\define \T {\Cal T}

\define\C{\Bbb C}

\def\today{\ifcase\month\or
       January\or February\or March\or April\or May\or June\or
       July\or August\or September\or October\or November\or December\fi
       \space\number\day, \number\year}

\NoBlackBoxes
\nologo

\loadeusm
\def\script#1{{\fam\eusmfam\relax#1}}
\input amssym
\nologo
\NoBlackBoxes

\topmatter

\title   Fourier transforms of positive definite kernels  and the Riemann $\xi$-Function 
\endtitle
\rightheadtext{ Fourier transforms and the Riemann $\xi$-Function}
\author George Csordas\endauthor

\address Department of Mathematics, University of Hawaii, Honolulu, HI
96822\endaddress
\email george\@math.hawaii.edu \endemail
\abstract 
The  purpose of this paper is to investigate the distribution of zeros of entire functions which can be
 represented as the Fourier transforms of certain admissible kernels. The principal results bring to light the intimate connection between the Bochner-Khinchin-Mathias theory of positive definite kernels and the generalized  real Laguerre inequalities. The concavity and convexity properties of the Jacobi theta function play a prominent role throughout this work.  The paper concludes with several questions and open problems.
\endabstract

\dedicatory{{Dedicated to Dr. Robert D. Cushnie on the occasion of his $80^{th}$ birthday.}}\enddedicatory

\endtopmatter

\baselineskip=18pt plus1pt minus 1pt
\parskip=10pt

\document
\centerline{{\bf 1. Introduction}}
  Today, there are no known explicit 
necessary and sufficient conditions that  even a ``nice"  kernel (cf. Definition 1.2), $K(t)$, must satisfy in order 
that its Fourier transform 
$$
F(x):=\frac{1}{2}\int_{-\infty}^{\infty} K(t) e^{ixt}\,dt= \int^\infty_0 K(t)\cos(xt)dt \tag 1.1
$$
have only real zeros (cf.  \cite{49, p. 17} and \cite{50}).   The program of investigation promulgated  here is motivated, in part, by several  recent results (\cite{4, 12, 24, 32, 42--45}) and our understanding that it is desirable to discover properties of the kernel, $K$, which (hopefully) will lead to information about the distribution of zeros of the entire function $F$.  The main leitmotif of this note pertains to  certain inequalities, known as the {\it generalized Laguerre inequalities} (Section 2), which play a pivotal role in the study of functions in the \LP class (cf.  Definition 1.1). Notwithstanding the extensive  research in this area and the impressive results dealing with the Riemann $\xi$-function, it is curious that to date so little progress has been made in proving some of the simplest Laguerre inequalities that $F$ must satisfy in order that it possess only real zeros (cf. Open Problem 4.7).

An outline of this work is as follows. In the remainder of this introduction, we recall some pertinent definitions and nomenclature that will  be used in the sequel.  In Section 2, we review several classical \underbar{and} new results involving the Laguerre and the generalized real Laguerre inequalities (Theorem 2.4) and prove two important, albeit elementary, results (Propositions 2.2 and 2.3) which adumbrate some of the applications in Section 4.
 With the aid of  the classical theorems of  S. Bochner \cite{1}, A. Khinchin \cite{33}, and M. Mathias \cite{41}, we establish the positive definite character of certain  canonical kernels which lead to some new classes of characteristic functions (Section 3). By extending the work of J. L. W. V. Jensen \cite{31} and G. P\'olya \cite{49}, our main results in Section 3 (cf. Theorems 3.5--3.7) establish precise relationships between certain positive definite kernels and the generalized real Laguerre inequalities. Concavity plays a prominent role throughout this paper and it is the {\it sine qua non} for analyzing the Jacobi theta function and related kernels. In Section 4, we apply the foregoing results and derive new necessary and sufficient conditions for the Fourier transform of the Jacobi theta function, the Riemann $\xi$-function, to belong to the Laguerre-P\'olya class. The paper ends with several (6) open problems (Section 4).

 In the present investigation, we will  adopt the following notation and nomenclature associated with real entire functions whose zeros lie in a strip. Let $S(\tau)$ denote the closed strip
of width $2\tau$, $\tau\ge 0$, in the complex plane, $\C$, symmetric  
about the real axis:
$$
S(\tau)=\{z\in \Bbb {C} \mid |\im (z) | \le \tau\}. \tag 1.2
$$

\proclaim{Definition 1.1}  We say that a real entire function $f$
belongs to the class $\s$, if $f$ can be expressed in the form
$$
f(z)=C e^{-az^2+bz}z^m\prod_{k=1}^\omega (1-z/z_k)e^{z/z_k},\qquad(0\le \omega\le \infty),\tag 1.3
$$
where $a\ge 0$, $b\in \Bbb R$, $z_k \in {S(\tau)}\setminus\{0\}$,
$\sum_{k=1}^\infty 1/|z_k|^2 < \infty$. We  allow functions in $\s$ to have only
finitely many zeros by letting, as  usual, $z_k=\infty$ and  $0=1/z_k$,
$k\ge k_0$, so that the canonical product in (1.3) is a finite product. By convention, the empty product is one.
If  $f\in\s$, for
some $\tau\ge0$, and if $f$ has only real zeros (i.e., if $\tau=0$),
then $f$ is said to belong to the
{\it Laguerre-P\'olya class}, and we write $f\in \lp$.  In addition, we  
write
$f\in \lp^*$, if $f=pg$, where $g\in \lp$ and $p$
is a real polynomial.  Thus, $f\in \lp^*$ if and only if $f\in\s$,
for some $\tau\ge 0$, and $f$ has at most finitely many non-real zeros.
\endproclaim
The significance
of the class $\s$ in the theory of entire functions stems from the fact that
$f\in \s$ if and only if $f$ is the uniform limit on compact sets
of a sequence
of real polynomials having zeros only in the
strip $S(\tau)$ (cf. \cite{2, p. 202} and  \cite{35b,  pp\, 373--374}).  It follows from the  Gauss-Lucas Theorem (\cite{40, pp 8--22}, \cite{53, p. 71}) that
this class of polynomials is closed under differentiation, and thus so
is $\s$.  For various properties and
algebraic and transcendental characterizations of functions in  the Laguerre-P\'olya class,
we refer the reader to P\'olya and Schur (\cite{52, p. 100},  \cite{51}, 
\cite{46, Kapitel II} or \cite{38, Chapter VIII}).

In the sequel, we will confine our attention to  special kernels which we term admissible kernels and define as follows.

\definition{Definition 1.2}  A function $K : \Bbb R \lra \Bbb R$ is called an 
{\it admissible kernel}, if it satisfies the following properties:
(i)
$K(t)\in C^{\infty}(\Bbb R)$, (ii)
$K(t) > 0$ for $t \in \Bbb R$, (iii)
$K(t) = K(-t)$ for $t \in \Bbb R$, (iv)
$K^{\prime}(t) < 0$ for $t > 0$, and
(v) for some $\varepsilon > 0$ and $n = 0, 1, 2, \dots$,
$$
K^{(n)}(t) = O\left(\op{exp}\big(-|t|^{2 + \varepsilon}\big)\right) 
\text{ as } t \lra \infty\,.\tag1.4
$$
\enddefinition 
Thus, the assertions that  $F(x)$ (cf. 1.1)  is a real entire function readily follows if we assume that $K(t)$
is an admissible kernel. Moreover,  a calculation shows
(\cite{52, p.\, 269}) that
$F(x)$ is an entire function of order $\frac{2+\varepsilon}{1+\varepsilon} <2$.  Also, 
by the Riemann-Lebesgue Lemma   $F(x)\to 0$ as $|x|\to \infty$. Observe that if we omit the requirement that $K(t)$ is even (see, Definition 1.2 (iii)), then its transform, $F$, cannot have only real zeros. This claim follows from integrating by parts and invoking the Riemann-Lebesgue Lemma (cf. \cite{49}).
\bigskip
\centerline{{\bf 2. The Laguerre Inequalities}} 
One important property, shared by all functions in $\lp$, is
logarithmic concavity; that is, if $f(x) \in \lp$, then $f(x)^2\left(\log
f(x)\right)^{\prime\prime}\le 0$ for all $x\in \Bbb{R}$. In order to
 verify this claim, one need only to consider the derivative of the logarithmic derivative of
$f(x)\in\lp$ using the (Hadamard) factorization  (1.3), (see, for
example,
\cite{6, 7, 8}). The logarithmic concavity, in
conjunction with the closure property of
$\lp$ under differentiation, implies that if
$f\in \lp$, then $f$ satisfies the following inequalities,
known as the {\it Laguerre inequalities}, (\cite{11--13, 15, 18, 22, 25, 47})
$$
L_{1,p}(x; f):=(f^{(p)}(x))^{2} -  f^{(p-1)}(x)
f^{(p+1)}(x)\geq 0,  \quad p=1, 2, 3, \dots,\quad \text{for
all}\quad x\in
\Bbb R. \tag 2.1
$$
For the sake of simplicity of notation, we  set  $L_{1,1}(x;f):=L_1(x;f):=L_1(x)$. In the sequel, we will be primarily concerned with the case when $p=1$ in (2.1); that is, $L_1(x)$. The reason for the subscript ``1" will become clear when we consider the generalized  real Laguerre inequalities (see Theorem 2.4). We remark that one of
the simplest manifestations of the existence of a non-real
zero of an entire function $f$,  occurs  when $f$ possesses a positive local minimum or a negative local maximum.  It is this observation that motivates us to consider the Laguerre inequalities.  We emphasize here that the Laguerre inequalities are only  {\it necessary}
conditions and, in general,  are not sufficient for an entire function to have only real
zeros. Indeed,  $f(x):=e^{-x^2}(1+x^2)\notin \lp$,
while a calculation shows that 
$L_1(x)= 2e^{-2x^2}x^2(3+x^2)\ge 0$ for all $x\in \Bbb R.$
\medskip\noindent
{\bf Remark 2.1.} To illustrate by an example the spirit of the type of research program  we are advocating here, consider again an admissible kernel $K(t)$ and its Fourier  cosine transform $F(x)$. Then via the change of variables, $u=-x^2$, we obtain the entire function 
$$
 F_c(u):=\sum_{k=1}^{\infty}\frac{k!\beta_k}{(2k)!}\frac{u^k}{k!}:=\int_0^{\infty}K(t)\cosh(t\sqrt{u})\,dt,\quad\text{where}\quad \beta_k:=\int_0^{\infty}K(t)t^{2k}\,dt, \quad k=0,1,2,\dots.
$$
Now set $\gamma_k:= \frac{k!\beta_k}{(2k)!}$ for $ k=0,1,2,\cdots$. If $\log(K(\sqrt t ))$ is strictly concave for all $t >0$, then we can infer that the Taylor coefficients of $F_c(x)$ satisfy the {\it Tur\'an inequalities}; that is, $L_{1,p}(0; F_c):=(F_c^{(p)}(0))^{2} -  F_c^{(p-1)}(0)F^{(p+1)}(0)=\gamma_p^2-\gamma_{p-1}\gamma_{p+1}\ge 0$,  for $p=1,2,3, \dots$ (see, for example, \cite{8, 14, 16, 21}).
Once again, the Tur\'an inequalities are only necessary conditions for $F_c$ (and whence for $F$) to belong to the Laguerre-P\'olya class.
\medskip
Our next proposition asserts that if a real entire function $f\in\s$,
$\tau=1$, has only real zeros in a vertical strip $A\le \re z \le B$,
$B-A>2$,  then
$L_1(x)\ge 0$ for
$x\in [A+1, B-1]:=I$. Thus, on the interval $I$,  $f$ cannot have a positive
local minimum or a negative local maximum. 
\proclaim{Proposition 2.2}  {\rm (\cite{22})}
  Let $f\in \s$, where $\tau=1$ and suppose that
$f(0)\not=0$. Let
$\{x_k\}_{k=1}^{\infty}$ denote the real zeros and let
$z_j=\alpha_j+i\beta_j$, $j=1,2, \dots,\omega$, $1\le\omega\le\infty$, denote
the non-real zeros of $f$. If there is an interval $[A, B]$,
with $B-A>2$, such that $\alpha_j \notin [A,B]$ for all $j\ge 1$, then
$$
L_1(x)\ge 0\quad \text{for all}\quad x\in [A+1, B-1].\tag 2.2
$$
\endproclaim
\demo{Proof} Using the product representation (1.3), logarithmic
differentiation yields
$$
L_1(x)=(f(x))^2\left\{2\alpha+\sum_{k=1}^{\infty}\frac{1}{(x-x_k)^2}
+2\sum_{j=1}^{\omega}\frac{(x-\alpha_j)^2-\beta_j^2}
{[(x-\alpha_j)^2+\beta_j^2]^2}\right\}. \tag 2.3
$$
Since $(x-\alpha_j)^2-\beta_j^2>0$ for any  $x\in [A+1, B-1]$, (2.3) gives
the desired result (2.2).
\qed\enddemo
\proclaim{Proposition 2.3}  {\rm (\cite{22})} Let $g(x)$ be a real entire function and define
$$
f(x):=\left ((x-\alpha)^2+\beta^2\right)^m\,g(x)\quad (\alpha\in \Bbb R, \beta >0, m\in \Bbb N), \tag 2.4
$$
so that $\alpha \pm i\beta$ are two non-real zeros of order $m$ of $f$. If $g(\alpha)\not= 0$, then
$$
L_1(\alpha; f)=-2m\beta^{4m-2}( g(\alpha))^2+\beta^{4m}L_1(\alpha; g).\tag 2.5
$$
Thus,  there exits $M>0$ sufficiently small such that 
$$
L_1(\alpha; f)< 0\quad \text{for all} \quad 0 < \beta <M.  \tag 2.6
$$
\endproclaim
\demo{Proof} Since $f(\alpha)=\beta^mg(\alpha)$, a straightforward calculation, using logarithmic differentiation,
yields  (2.5) and whence the desired result (2.6) follows.
\qed\enddemo
\noindent
 A heuristic description of Proposition 2.2 is as follows. A  conjugate pair of non-real zeros $\alpha\pm i\beta$ of $f(x)$, when $\beta >0$ is sufficiently small, forces $L_1(\alpha; f)$ to be negative.

We   consider next the so-called  {\it generalized real Laguerre inequalities} (see, for example, \cite{20, 25}) that are {\it both necessary and sufficient} for membership in the \LP class.
\noindent
\proclaim{Theorem 2.4} {\rm(The Generalized Real Laguerre Inequalities \cite{20, Theorem 2.9})}
Let  $f$ denote a real entire function,  $f \not  \equiv 0$. For $n\in \Bbb N_0:=\Bbb N\cup\{0\}$ and $x\in \Bbb R$, set
$$
L_n(x):=L_{n,1}(x; f):=\sum_{j=0}^{2n} \frac{(-1)^{j+n}}{(2n)!} \binom{2n}{j}
f^{(j)}(x)f^{(2n-j)}(x). \tag 2.7
$$
$$
\text{If}\quad f(x)\in \lp,\quad\text{then}\quad L_n(x)\ge 0\quad  \text{for all}\quad n\in\Bbb N_0\quad\text{and for all } x\in\Bbb R. \tag 2.8
$$
Conversely, suppose  that  
$$
f(x)=e^{-a x^2}g(x),\quad  a \ge 0, \text{ where the  genus of }g(x) \text{ is } 0 \text{ or } 1. \tag 2.9
$$
$$ 
\text{If} \quad L_n(x)\ge 0\quad  \text{for all}\quad n\in\Bbb N_0\quad\text{and for all } x\in\Bbb R, \quad \text{then}\quad f(x)\in \lp. \tag 2.10
$$
\endproclaim
\noindent
{\bf Remarks 2.5.} Observe that $L_0(x)=f(x)^2$ and  to justify the appellation ``generalized Laguerre expression", note that $L_1(x)=f^{\prime}(x)^2-f(x)f^{\prime\prime}(x)$. In  addition, we remark that if the real entire function $f(x)$ satisfies the generalized real Laguerre inequalities, $L_n(x)\ge 0$ ($n\in\Bbb N_0,\, x\in\Bbb R$), then $f(x)$ has only real zeros (cf. \cite{20, p.\, 343}). For the sake of completeness, we mention here the following representation of $|f(x+iy)|^2$ which  can be derived by a direct calculation (see, for example, \cite{20}, \cite{47},  \cite{49} or  by using a recursion relation \cite{6}): 
$$
|f(x+iy)|^2=f(x+iy) f(x-iy)=\sum_{n=0}^{\infty}L_n(x)y^{2n}, \qquad (x,\,y \in \Bbb R), \tag 2.11
$$
where $L_n(x)$ is defined in (2.7).

\noindent
{\bf Remarks 2.6.} The action of the non-linear operators $\{L_n\}_{n=0}^{\infty}$ taking a real entire function $f(x)$ to $L_n(x; f):=L_n(x)$ is given implicitly by equation (2.11). We mention here, parenthetically, a couple facts about these operators. It is known that the operators $L_n$ satisfy  a simple recursive relation \cite{6, Theorem 2.1} and that $L_n(x)$ is also a real entire function \cite{6, Remark 2.4}. Interesting generalizations of these operators are given by K. Dilcher and K. B. Stolarsky \cite{26} and D. A. Cardon \cite{3} (see also Section 3). 
Recently, A. Vishnyakova and the author \cite{25}  have shown that the various sufficient conditions for a real entire function, $f(x)$, to belong to the \LP class, expressed in terms of Laguerre-type inequalities, do not require the {\it a priori} assumptions  about the order and type of $f(x)$. Thus, for instance,  implication (2.10) remains valid if we omit assumption (2.9). In light of the results in \cite{25}, we can state the {\it complex Laguerre inequalities} as follows. Suppose $f$, $f\not\equiv 0$, is a real entire function. Once again we do not stipulate conditions on the order and type of $f$ \cite{25}. Then $f\in \lp$ if and only if  
$$
|f^{\prime}(z)|^2\ge \Re\left(f(z)
\overline{f^{\prime\prime}(z)}\right)\quad\text{for all}\quad z\in
\Bbb C. \tag 2.12
$$
It may be of interest to note  that the complex Laguerre expression can be also formulated in terms of  two real Laguerre-type expressions \cite{25}. Indeed, if $f(x+iy)=U(x,y)+iV(x,y)$ is a real entire function, then a calculation shows that for all $z=x+iy\in \Bbb C$,
$$
\frac{1}{2}\frac{\partial^2}{\partial y^2}|f(x+iy)|^2=|f^{\prime}(z)|^2- \Re\left(f(z)
\overline{f^{\prime\prime}(z)}\right)=U_x^2-UU_{xx}+V_x^2-VV_{xx}.\tag 2.13
$$
\bigskip
\centerline{{\bf 3. Positive Definite Functions and the Laguerre Inequalities}} 

 We  mention  at the outset that it was M. Mathias \cite{41} who in 1923, motivated by the results of C. Carath\'eodory and O. Toeplitz (cf. \cite{55, p.\,412}),  first defined and studied the properties of \pd  functions. In this section, after reviewing some definitions, we will succinctly summarize a couple of classical results due to M. Mathias \cite{41},  S. Bochner \cite{1}, A. Khinchin \cite{33} and G. P\'olya \cite{48}. Parenthetically we note that there are many excellent treatises in the literature dealing with positive definiteness and here we merely cite F. Lukacs \cite{39}, T. Kawata \cite{34}, M. Mathias \cite{41} and J. Stewart \cite{55}, together with the original works of S. Bochner \cite{1}, A. Khinchin \cite{33} and G. P\'olya \cite{48}. The interested reader will find 125 additional references in J. Stewart's outstanding survey article \cite{55}. In the second part of this short section, our goal is  to bring to  light  the connection between positive definiteness and the Laguerre inequalities.
\definition{Definition 3.0}\rm{(\cite{34, p.\, 377})} A continuous function $\varphi: \Bbb R \to \Bbb R$ is said to be {\it \pd } (or more precisely {\it non-negative definite}), if 
$$
\int_{-\infty}^{\infty} \int_{-\infty}^{\infty}\varphi (t-s) \rho(t)\overline{\rho(s)}\, dt\,ds\ge 0, \tag 3.1
$$
where $\rho: \Bbb R \to \Bbb C$ is  any measurable function with compact support.
\enddefinition

\noindent
An equivalent definition of positive definiteness \rm{(\cite{34, p.\, 377})} is the following discrete formulation. A continuous function $\varphi$ is \pd if the Hermitian form
$$
\sum_{j=1}^n\sum_{k=1}^n \varphi(x_j-x_k)\rho_j\overline{\rho_k} \ge 0 \quad\text{for every  }x_1,\dots,x_n\in\Bbb R \text{   and   }  \rho_1,\dots,\rho_n \in \Bbb C. \tag 3.2
$$
\noindent
By way illustration, we note that $\varphi(x)=\cos x$ is positive definite, since
 $$
 \sum_{j=1}^n\sum_{k=1}^n \cos (x_j-x_k)\rho_j\overline{\rho_k}=\left |\sum_{j=1}^n \rho_j\cos x_j\right |^2+\left |\sum_{k=1}^n \rho_k\sin x_k\right |^2\ge 0 \quad (x_1,\dots,x_n\in\Bbb R,  \quad \rho_1,\dots,\rho_n \in \Bbb C).
$$
Similarly, it is easy to check that $e^{itx}$, ($t\in \Bbb R$), is positive definite; while it is not so straightforward to verify that the functions $e^{-|x|}, e^{-x^2}$ and $\frac{1}{1+x^2}$ are positive definite.
\noindent
For the sake of clarity, we define one more term. By a {\it distribution function} we  shall mean a non-decreasing function $V(x)$ such that $V(-\infty)=0$ and $V(+\infty)=1$. The Fourier-Stieltjes transform of $V$, 
$$
f(t)= \int_{-\infty}^{\infty}e^{itx} dV(x) \qquad (-\infty <t <\infty), \tag 3.3
$$
is called the {\it characteristic function}  corresponding to the given distribution function $V$.

In 1932, S. Bochner proved the following celebrated theorem that bears his name.
\proclaim{Theorem 3.1}  {\rm (\cite{1}, \cite{39, p.\, 71})} A continuous function, $f(t)$, with $f(0)=1$, is a characteristic function if and only if $f(t)$ is positive definite.
\endproclaim
\noindent
We remark that since  $e^{itx}$, ($t\in \Bbb R$) is positive definite, it is easy to show that a characteristic function is positive definite. The converse implication is the difficult part of Theorem 3.1 (see, for example, T. Kawata \cite{34,  p.\, 377} or E. Lukacs \cite{39, p.\,71}). For our purposes the following version  of the Khinchin's criterion \cite{33} for a characteristic function will suffice (see also E. Lukacs \cite{39, Theorems 4.2.4 and 4.2.5}).
\proclaim{Theorem 3.2}  {\rm (\cite{34, p.\,387})} A  function of the form
$$
f(t)=\frac{1}{c}\int_{-\infty}^{\infty}\varphi(x+t)\overline{\varphi(x)}\,dx, \tag 3.4
$$
where $\varphi(x)$ is any function in $L^2(\Bbb R)$ with $||\varphi ||_2=c >0$, or the local uniform limit of such functions, is a characteristic function. The converse is also true.
\endproclaim
\noindent
Theorem 3.2 implies Mathias's result \cite{41, Satz 15} which may be stated as follows. 
 If $\varphi\in L^2(\Bbb R)$, then the function
 $$
f(t)=\int_{-\infty}^{\infty}\varphi(s+t)\overline{\varphi(s-t)}\,ds, \tag 3.5
$$
is positive definite. We remark that if $\varphi: \Bbb R \to \Bbb R$ is an admissible kernel, then $\varphi$ is a bounded integrable function. Moreover, it is not difficult to demonstrate that $\varphi$ satisfies the conditions of Fourier's inversion theorem (cf. \cite{56, Pringsheim's theorem, p.\, 16}). Thus, with the terminology adopted here we can express Mathias's main theorem (cf. \cite{41, Hauptsatz, p.\, 108} or \cite{55, p.\, 412}) in the following form.
\proclaim{Theorem 3.3} Let $\varphi$ be an admissible kernel and let
$$
f(t):= \int_{-\infty}^{\infty}\varphi(x)\cos(xt)\,dx.\tag 3.6
$$
Then $\varphi$ is \pd if and only if $f(t)\ge 0$ for all $t\in \Bbb R$.
\endproclaim
\noindent
The above necessary and sufficient conditions for a characteristic function are, in general, not readily applicable in order to determine whether a given function is a characteristic function. There is, however, a beautiful and simple criterion due to P\'olya \cite{48} (see also Lukacs \cite{39,   p.\, 85}).
\proclaim{Theorem 3.4}  {\rm (P\'olya's criterion)} Suppose that $f: \Bbb R \to  \Bbb R$ is continuous and satisfies the following conditions:  (i) $f(0)=1$, (ii) $f(-t)=f(t)$, (iii) $f$ is convex for $t>0$ and (iv) $\lim_{t\to\infty}f(t)=0$. Then $f(t)$ is the characteristic function of an absolutely continuous distribution function $V(x)$.
\endproclaim
\noindent
Thus, P\'olya's criterion provides a sufficient condition for a continuous function $f: \Bbb R \to  \Bbb R$  to be a characteristic function. There is however a caveat in order:  our admissible kernels do not satisfy the convexity hypothesis of Theorem 3.4.
Preliminaries aside, we will now  relate positive definiteness to the  various Laguerre-type inequalities presented in Section 2. Our first result in this direction shows that if $\varphi(t)$ is an admissible kernel such that $\log \varphi(t)$ is  strictly concave (i.e., $d^2/dt^2 \log \varphi(t) < 0$ for $t> 0$), then for each $n\in \Bbb N \cup {0}$ we can  associate  with $\varphi(t)$ a (canonical) kernel $K_n$ which is also an admissible kernel.
\proclaim{Theorem 3.5} Let $\varphi(t)$ be an admissible kernel. If $\log \varphi(t)$ is  strictly concave  for $t> 0$, then for each non-negative integer $n$, the associated kernel
$$
K_n(t):=\int_{-\infty}^{\infty}\varphi(s+t)\varphi(s-t)s^{2n}\,ds \qquad (n=0,1,2,\dots), \tag 3.7
$$
is also an admissible kernel.
\endproclaim
\demo{Proof} Fix a non-negative integer $n$. Consulting Definition 1.2, we readily deduce  that $K_n(t)$ satisfies the properties (i), (ii), (iii) and (v) of Definition 1.2. Thus, it remains to show that  $K_n^{\prime}(t) < 0$ for $ t>0$. Invoking Leibniz's rule to justify the differentiation under the integral, we have
$$
K_n^{\prime}(t)=2\int_{0}^{\infty}\left[\varphi^{\prime}(t+s)\varphi(t-s)+\varphi(t+s)\varphi^{\prime}(t-s)\right]s^{2n}\,ds. \tag 3.8
$$
Next, we fix $t>0$ and consider the intervals of integration $I_1:=(0,t)$ and $I_2:=(t,\infty)$. Since $\log \varphi(t)$ is  strictly concave for $t>0$, ${\ds \frac{\varphi^{\prime}(t)}{\varphi(t)}}$ is strictly decreasing for $t>0$. Hence, for $s>0$, we claim  that
$$
\frac{\varphi^{\prime}(t+s)}{\varphi(t+s)}< - \frac{\varphi^{\prime}(t-s)}{\varphi(t-s)}. \tag3.9
$$
If $s\in I_1$, then $0< s<t$ and $-\varphi^{\prime}(t-s)>0$. Since $\varphi^{\prime}(t+s)<0$, we see that (3.9) holds. On the other hand, if  $s\in I_2$, then $t-s < 0$. Since $\varphi(t)$ is an even function, $\varphi(t-s)=\varphi(s-t)$. Also, $0 < s-t < s+t$, and thus, (3.9) holds, since
$$
 \frac{\varphi^{\prime}(t+s)}{\varphi(t+s)}< \frac{\varphi^{\prime}(s-t)}{\varphi(s-t)}=-\frac{\varphi^{\prime}(t-s)}{\varphi(t-s)}. 
 $$
 \qed\enddemo
\noindent
Following P\'olya's work involving Jensen's {\it Nachlass} (\cite{49, pp 278--308}), we next establish an important relationship between  a given strictly logarithmically concave admissible kernel and the associated admissible kernel $K_n(t)$ defined in (3.7).
\proclaim{Lemma 3.6} If $\varphi(t)$ is a strictly logarithmically concave admissible kernel for $t>0$, then
$$
\int_{-\infty}^{\infty}\int_{-\infty}^{\infty}\varphi(t)\varphi(s)e^{ix(s+t)}(s-t)^{2n}\,dt\,ds=2\cdot 2^{2n}\int_{-\infty}^{\infty}K_n(v)\cos(2xv)\,dv, \qquad (n=0,1,2,\dots)\tag 3.10
$$
where $K_n$ is the associated admissible kernel defined by (3.7).
\endproclaim
\demo{Proof} (A sketch.)  Consider the entire function
$$
F(x):=\int_{-\infty}^{\infty}e^{itx} \varphi(t)\,dt. \tag 3.11
$$
Then, since  both $\varphi(t)$ and $K_n(t)$ are admissible kernels (Theorem 3.5), the following calculations are valid: 
$$
\align
|F(x+iy)|^2&= F(x+iy)F(x-iy)\\
&= \int_{-\infty}^{\infty}\int_{-\infty}^{\infty}\varphi(t)\varphi(s)e^{ix(s+t)}e^{-(s-t)y}\,dt\,ds\\
&=\sum_{n=0}^{\infty}\frac{y^{2n}}{(2n)!}\int_{-\infty}^{\infty}\int_{-\infty}^{\infty}\varphi(t)\varphi(s)e^{ix(s+t)}(s-t)^{2n}\,dt\,ds.\tag 3.12
\endalign
$$
Next, we  use (i) Euler's formula $e^{ix(s+t)}=\cos (x(s+t))+i\sin(x(s+t))$, (ii) the fact that the odd functions integrate to zero and (iii) the absolute value of the Jacobian of the transformation, $s\to u+v$ and $t\to u-v$, is 2. Accordingly, we obtain
$$
\align
&\int_{-\infty}^{\infty}\int_{-\infty}^{\infty}\varphi(t)\varphi(s)e^{ix(s+t)}(s-t)^{2n}\,dt\,ds\\
&=\int_{-\infty}^{\infty}\int_{-\infty}^{\infty}\varphi(t)\varphi(s)\cos(x(s+t))(s-t)^{2n}\,dt\,ds\\
&=\int_{-\infty}^{\infty}\int_{-\infty}^{\infty}\varphi(t)\varphi(s)\cos(x(s-t))(s+t)^{2n}\,dt\,ds\\
&=2\cdot 2^{2n}\int_{-\infty}^{\infty}\int_{-\infty}^{\infty}\varphi(u+v)\varphi(u-v)\cos(2xv)\,(u)^{2n}\,du\,dv\\
&=2\cdot 2^{2n}\int_{-\infty}^{\infty}K_n(v)\cos(2xv)\,dv. \tag 3.13
\endalign
$$
\qed\enddemo
\proclaim{Theorem 3.7} Let $\varphi(t)$ be a strictly logarithmically concave (for $t>0$) admissible kernel and let
$K_n$ ($n=0,1,2,\dots$) denote the associated admissible kernel defined by (3.7). Let $F(x):=\int_{-\infty}^{\infty}e^{itx} \varphi(t)\,dt.$  Then,
$$
L_n(x):=L_n(x; F):=\frac{2\cdot2^{2n}}{(2n)!}\int_{-\infty}^{\infty}K_n(t)\cos(2xt)\,dt, \qquad (n=0,1,2,\dots),\tag 3.14
$$
where $L_n(x)$ is the generalized real Laguerre expression (cf. (2.7) of Theorem 2.4) for the entire function $F$. Moreover, $F\in \lp$ if and only if $K_n$ is a \pd kernel for all $n=0,1,2,\dots$.
\endproclaim
\demo{Proof} We recall from Section 2 (see (2.11) of Remark 2.5) that the Taylor coefficient of $y^{2n}$ in the expansion $|F(x+iy)|^2$ is precisely $L_n(x)$. Hence, (3.14) follows from (3.12) and (3.13). Next, by Theorem 2.4, $F\in \lp$ if and only if $L_n(x)\ge 0$ for $n=0,1,2,\dots$ and for all $x\in \Bbb R$. By Theorem 3.3 and (3.14), $L_n(x)\ge 0$ for $n=0,1,2,\dots,$ and for all $x\in \Bbb R$, if and only if  the kernel $K_n$ ($n=0,1,2,\dots$) is positive definite.
\qed\enddemo
\noindent
{\bf Remarks 3.8.} (a) In \cite{41, Satz 15}, Mathias has proved that the kernel $K_0(t)=\int_{-\infty}^{\infty}\varphi(s+t)\varphi(s-t)\,ds$ is positive definite. This is clear in our setting, since $L_0(t)=|F(t)|^2$. The case when $n=1$; that is, 
$$
K_1(t)=\int_{-\infty}^{\infty}\varphi(s+t)\varphi(s-t)\, s^2\,ds \quad\text{and}\quad  L_1(x)=(F^{\prime}(x))^2-F(x)F^{\prime\prime}(x)=4\int_{-\infty}^{\infty}K_1(t)\cos(2xt)\,dt,
$$
appears to be much more difficult. (b) The desideratum to characterize Fourier transforms in the Laguerre-P\'olya class, in terms of the indicated kernels,  is  achieved by Theorem 3.7.  However, the elusive nature of positive definiteness certainly remains as an issue. (c) It may be noteworthy to remark  that in conjunction with the Bochner and Khinchin results (cf. Theorems 3.1 and 3.2), our Theorem 3.7 gives rise to new families of characteristic functions when the kernels $K_n$ associated with functions $F\in\lp$ are appropriately normalized.

\noindent
{\bf Open Problem 3.9.} Characterize the logarithmically concave admissible kernels $\varphi(t)$ such that the  associated admissible kernel $K_1(t)$ is positive definite.

Striving for simplicity, we propose here another, direct, approach for showing that $\int_{0}^{\infty}K_1(t)\cos(xt)\,dt\ge 0$ for all $x\in\Bbb R$.
\proclaim{Proposition 3.10} Let $F(x):=\int_{0}^{\infty} \varphi(t)\cos xt\,dt$ and set $K_1(t):=\int_{0}^{\infty}\varphi(s+t)\varphi(s-t)\, s^2\,ds$. Let $\overline{G}(t):=\int_{t}^{\infty}K_1(u)\,du$ and  $A:=\overline{G}(0)$. Then $K_1(t)$ is \pd if and only if
$$
\int_0^{\infty}\overline{G}(t)\sin xt\, dt\le \frac{A}{x}\qquad \text{for all}\quad x \not = 0. \tag 3.15
$$
\endproclaim
\noindent
{\bf Remark 3.11.} Before we prove Proposition 3.10, we recall that P\'olya's argument \cite{48} shows that in general, the non-negativity of the Fourier sine transform is easier to demonstrate than that of the of the Fourier cosine transform. Indeed, consider the function  $\overline{G}(t)$ defined in Proposition 3.10. Then for each fixed $x>0$,
$$
\align
I(x):&=\int_0^{\infty}\overline{G}(t)\sin xt\,dt=\sum_{k=0}^{\infty}\int_{\pi k/x}^{\pi (k+1)/x}\overline{G}(t)\sin xt\,dt\qquad \left(t=s+\frac{\pi k}{x}\right)\\
&=\sum_{k=0}^{\infty}\int_{0}^{\pi /x}\overline{G}\left(s+\frac{\pi k}{x}\right)\sin (xs+\pi k)\,ds\\
&=\sum_{k=0}^{\infty}(-1)^k\int_{0}^{\pi /x}\overline{G}\left(s+\frac{\pi k}{x}\right)\sin (xs)\,ds.
\endalign
$$
Since $\overline{G}(s)>0$,  $\overline{G}^{\prime}(s)<0$ ($s>0$), and $\overline{G}(s)\to 0$ as $s\to \infty$, it follows from the alternating series test that $I(x)>0$ for $x>0$.
\demo{Proof of Proposition 3.10} Integration by parts yields,
$$
\align
\int_{0}^{\infty}K_1(t)\cos(xt)\,dt&=\int_{0}^{\infty}K_1(u)\,du-x\int_{0}^{\infty}\left(\int_{t}^{\infty}K_1(u)\,du\right)\sin xt \,dt\\
& = A-\int_0^{\infty}\overline{G}(t)\sin xt\, dt\\
\endalign
$$
and whence inequality (3.15) follows if and only if $K_1(t)$ is positive definite.
\qed\enddemo 
We conclude this section with a concrete example which demonstrates that,  if $K_n$ is positive definite, then 
in general, $K_{n+1}$ need not be positive definite. There are several ways we can illustrate this fact. The kernel we will use is a  Gaussian, $e^{-t^2}$, times a polynomial  and therefore it will not satisfy condition (v) of Definition 1.2. Nevertheless,  our choice facilitates the exact evaluation of the required integrals.  The calculations are sufficiently involved, albeit elementary, to warrant the use of a computer. 

\noindent
{\bf Example 3.12.}  Let $\varphi(t):=e^{-t^2}(15+t^2+t^4)$.  Then it is easy to confirm  that $\varphi(t)$ satisfies conditions (i)--(iv) (but not (v)) of Definition 1.2. In addition, $\log(\varphi(t))$ is strictly concave for $t>0$. In the subsequent calculations, we will denote by $c_j$, $j\ge 1$, a positive constant whose exact value is irrelevant. Then
$F(x)=\int_{-\infty}^{\infty} \varphi(t)\cos xt\,dt= c_1e^{-x^2/4} (260-16x^2+x^4)$. Since $F(x) >0$, $F$ has 4 non-real zeros (i.e., $F\notin\lp$) and whence by Theorem 3.7 at least one of the kernels $K_n$ (cf. (3.7)) fails to be positive definite. Since $\int_{-\infty}^{\infty}K_1(t)\cos 2xt\,dt=c_2 e^{-x^2/2}(84240 - 13536 x^2 + 712 x^4 - 24 x^6 + x^8)> 0$ for all $x\in \Bbb R$, $K_1$ is (strictly) positive definite. On the other hand, $\int_{-\infty}^{\infty}K_2(t)\cos 2xt\,dt=c_3e^{-x^2/2}(107088 - 18496 x^2 + 696 x^4 - 16 x^6 + x^8)$ has 4 simple real zeros and consequently
$K_2$ is {\it not} positive definite.
\bigskip
\centerline{{\bf 4. Scholia: the Jacobi Theta Function  and the Riemann $\xi$-function}}

  The purpose of this section is three-fold: (i) to investigate the properties of the Jacobi theta  function (cf. (4.2)) and related kernels, (ii) apply the results of Section 3  (Theorem 3.5 and Theorem 3.7) and provide new necessary and sufficient conditions for  $H(x):=\xi(x/2)/8\in \lp$ (cf. (4.1)), and (iii) formulate some open problems involving kernels associated with the Jacobi theta  function.

By way of background information,  we  commence with Riemann's definition of his $\xi$-function (\cite{49, p.\,10}); that is, 
$$
\xi(iz):=\frac{1}{2}\left(z^2-\frac{1}{4}\right)\pi^{-z/2-1/4}\Gamma \left(\frac{z}{2}+\frac{1}{4}\right)\zeta\left(z+\frac{1}{2}\right),
$$
Then it is known  (\cite{49 p.\,11}), \cite{57, p.\,255}  or \cite{52, p.\, 286})
that $\xi(x)$ admits the integral representation of the form 
$$
H(x): = \frac18\xi\left(\frac x2\right): = \int^\infty_0 \Phi(t)\cos(xt)dt\,,
\tag4.1
$$
where  the {\it Jacobi theta function}, (without the usual factor 4) is defined as
$$
\Phi(t): = \sum^\infty_{n=1}\pi n^2\big(2\pi n^2e^{4t} - 3\big)\op{exp}\big(5t - 
\pi n^2e^{4t}\big)\,. \tag 4.2
$$
The Riemann Hypothesis is equivalent to the statement that all the  
zeros of
$H(x)$ are real (cf. \cite{57, p. 255}). We also recall that $H(x)$ is  
an entire
function of order one \big(\cite{57, p. 16}\big) of maximal type  (cf.
\cite{19, Appendix A}).  Thus, with the above nomenclature (cf. Section 1) the Riemann
Hypothesis is true if and only if
$H\in\lp$. It is also known  (\cite{30, p.\, 7} ) that all the zeros of
$H$ lie in the {\it interior} of the strip $S(1)$, so that $H(x)\in  
\s$, with $\tau=1$ and  that $H(x)$  has an infinite  number of real zeros \cite{57, p. 256}.
Before we begin with  a synopsis of results, we emphasize that  the {\it raison d'\^etre} for investigating the kernel $\Phi$ is that there is an intimate connection (the precise meaning of which is yet unknown) between the properties of $\Phi$  and the distribution of the zeros of its Fourier transform $H(x)$ (cf. (4.1)).  
\proclaim{Theorem 4.1}  {\rm (\cite{16, Theorem A})}  Consider the 
function $\Phi$ of {\rm (4.2)} and set
$$
\Phi(t) = \sum^\infty_{n=1} a_n(t)\,,\quad \text{where}\quad a_n(t):= \pi n^2\big(2\pi n^2e^{4t} - 3\big)\op{exp}\big(5t - \pi n^2e^{4t}\big)\qquad (n = 1, 2, \dots)\,.
$$
Then, the following are valid:
\roster
\noindent \item"{(i)}"
for each $n \ge 1$, $a_n(t) > 0$ for all $t \ge 0$, so that $\Phi (t) > 0$ for 
all $t \ge 0;$
\item"{(ii)}"
$\Phi(z)$ is analytic in the strip $-\pi/8 < \op{Im} z < \pi/8;$
\item"{(iii)}"
$\Phi$ is an even function, so that $\Phi^{(2m + 1)}(0) = 0 \quad (m = 0, 1, 
\dots);$
\item"{(iv)}"
for any $\varepsilon > 0$, $\lim_{t\to\infty}\Phi^{(n)}(t)\op{exp}\big[(\pi - \varepsilon)e^{4t}\big] =0$;
\item"{(v)}"
$\Phi'(t) < 0$ for all $t > 0$\,.
\endroster
\endproclaim
\noindent The proofs of statements (i) -- (iv) can be found in G. P\'olya \cite{49}, whereas 
the proof of (v) is in A. Wintner \cite{58} (see also Spira \cite{54}).  \qquad$\square$

In order to indicate the significance of the next theorem, we consider the Taylor series  of $H(x)$ about the origin
$$
H(z) = \sum^\infty_{k=0}\frac{(-1)^ k b_k}{(2k)!}z^{2k}\,,\quad\text{where}\quad  b_k := \int^\infty_0t^{2k}\Phi(t) dt\qquad(k = 0, 1, 2, \dots).
\tag 4.3
$$
The change of variable, $z^2 = -x$ in (4.3), yields the entire function
$$
F(x) := \sum^\infty_{k=0}\frac{\gamma_k}{k!}x^k, \quad \text{where}\quad \gamma_k:=\frac{k!b_k}{(2k)!}>0\quad (k=0,1,2,\dots). \tag4.4
$$
Then it is easy to see that $F(x)$ is an entire function of order $\frac12$ 
and that the Riemann Hypothesis  is equivalent to the 
statement that all the zeros of $F(x)$ are real and negative.  Now it is 
known  (P\'olya and Schur \cite{51}) that a {\it necessary 
condition} for $F(x)$ to have only real zeros is that the moments $b_k$ 
(in (4.3)) satisfy the {\it Tur\'an inequalities}; that is, 
$$
 b_k^2 - \frac{2k - 1}{2k + 1} b_{k-1} b_{k+1} \ge 0 \quad\text{ or equivalently}\quad  T_k:=\gamma_k^2-\gamma_{k-1}\gamma_{k+1}\ge 0 \quad(k = 1, 2, 3, \dots)\,. \tag4.5
$$
These inequalities have been established (cf. \cite{16},  for $m\ge 2$, and 
\cite{21}) as a consequence of either one of the two concavity properties ((a) or (b)) 
of $\Phi$ stated in the following theorem. (For related interesting results see also D. K. Dimitrov and F. R. Lucas \cite{29} and D. K. Dimitrov \cite{27}, \cite{28}).
\proclaim{Theorem 4.2}  Let $\Phi$ be defined by {\rm (4.1)}.  Then $\Phi$  satisfies 
the following concavity properties. 

{\rm (a) (\cite{16, Proposition 2.1})} If
$$
K_{\Phi}(t) := \int^\infty_t\Phi\big(\sqrt u\big)du\qquad(t \ge 0)\,,
$$
then $\log K_{\Phi}(t)$ is strictly concave for $t > 0$; that is, $\frac{d^2}{dt^2}
\log K_{\Phi}(t) < 0 \quad \text{for} \quad t > 0.$

{\rm (b) (\cite{21, Theorem 2.1})}  The function $\log \Phi(\sqrt{t})$ is strictly 
concave for 
$t > 0$.\qquad$\square$
\endproclaim
\noindent
{\bf Remarks 4.3.} (a) A calculation shows that $\log \Phi(\sqrt{t})$ is strictly  concave for $t > 0$ if and only if $g(t):= 
 t\left[\big(\Phi^{\prime}(t)^2- \Phi(t)\Phi^{\prime\prime}(t)\right] +  \Phi(t)\Phi^{\prime}(t) > 0$ for $ t>0$. Since $ \Phi(t) >0$ and $\Phi^{\prime}(t)<0$ for $ t>0$, it is easy to check that the  inequality $g(t) >0$ is stronger than the assertion that $\log(\Phi(t))$ is  strictly concave  for $ t>0$. Indeed, the inequality $\Phi^{\prime}(t)^2- \Phi(t)\Phi^{\prime\prime}(t)>0$ does not imply, in general, the Tur\'an inequalities (4.5) (see, for example, \cite{5, Example 3.4}).
 
 \noindent
 (b) Since  $ \Phi(t) >0$ and $\Phi^{\prime}(t)<0$ for $ t>0$, we can also demonstrate that that the ``average value" of $H(x)$, the Fourier cosine transform of $\Phi$ (cf. (4.1)), is positive. Indeed, for $t>0$,
 $$
 \int_0^t H(u)\,du=\int_0^{\infty}\Phi(x)\left(\int_0^t\cos xu\, du\right)\,dx=\int_0^{\infty}\Phi(x)\frac{\sin xt}{x}\,dx>0,
 $$
where the last inequality can be established using the method of proof presented in Remark 3.11.
\noindent
We pause for a moment, and append here yet another convexity result involving $\Phi$.
\proclaim{Theorem 4.4}\rm{(\cite{11, pp 43--44})} The function $\Phi(\sqrt t)$ is strictly convex for $t>0$; (that is, $\frac{d^2}{dt^2}\Phi(\sqrt t)>0$ for $t>0$) and hence 
$$
\int_0^{\infty}\Phi(\sqrt t)\cos xt\,dt >0 \quad \text{ for all  } x\in \Bbb R.
 $$
\endproclaim
Having reviewed some of the salient properties of the Jacobi theta function, we are  now in position to apply the results of Section 3.
\proclaim{Theorem 4.5} The Jacobi theta function, $\Phi(t)$, is a strictly logarithmically concave admissible kernel. Moreover, the associated kernel
$$
K_n(t):=K_n(t; \Phi):=\int_{-\infty}^{\infty}\Phi(s+t)\Phi(s-t)s^{2n}\,ds \qquad (n=0,1,2,\dots), \tag 4.6
$$
is also an admissible kernel.
\endproclaim
\demo{Proof} By Theorem 4.1, $\Phi$ is an admissible kernel. Now, it follows from Theorem 4.2 and Remarks 4.3 (a) that $\log \Phi(t)$ is  strictly concave  for $t> 0$. Thus, by Theorem 3.5, for each non-negative integer $n$, the associated kernel $K_n(t):=K_n(t; \Phi)$ is also an admissible kernel.
\qed\enddemo
\noindent
Finally, with the aid of Lemma 3.6, Theorem 3.5 and Theorem 3.7, we obtain the following equivalent formulation of the Riemann Hypothesis.
\proclaim{Theorem 4.6} Let
$K_n:=K_n(t; \Phi)$ ($n=0,1,2,\dots$) denote the associated admissible kernel defined by (4.6). Let $H(x):=\int_{0}^{\infty} \Phi(t)\cos xt\, dt.$  Then, for $n=0,1,2,\dots$,
$$
L_n(x):=L_n(x; H):=\frac{2\cdot2^{2n}}{(2n)!}\int_{-\infty}^{\infty}K_n(t)\cos(2xt)\,dt, \tag 4.7
$$
where $L_n(x)$ is the generalized real Laguerre expression (cf. (2.7) of Theorem 2.4) for the entire function $H$. Moreover, $H\in \lp$ if and only if $K_n$ is a \pd admissible kernel for all $n=0,1,2,\dots$.
\endproclaim
At this juncture, we are obliged to expose our  ignorance and state the following tantalizing open problem.

\noindent
{\bf Open Problem 4.7.} (One of the simplest Laguerre inequalities for the Riemann $\xi$-function.) Let  $\Phi$ denote the Jacobi theta function and let $H(x):=\xi(x/2)/8=\int_{0}^{\infty} \Phi(t)\cos xt\, dt.$ Then, is it true that
$$
 L_1(x)=(H^{\prime}(x))^2-H(x)H^{\prime\prime}(x)\ge 0\quad \text{ for all}\quad x\in \Bbb R? \tag 4.8
$$
\noindent
{\bf Remark 4.8.} The verification of the special Laguerre inequality  (4.8) itself would be 
significant. If we could prove that $L_1(x) >0$ for all real $x$, then it would follow that all the real zeros of $H$ are {\it simple}. Of course, should inequality 
(4.8) fail to hold for some $x_0$, then the Riemann Hypothesis would be false.  Now it follows from the numerical results of van de Lune, te Riele, and Winter \cite{37} that the zeros of $H(x)$ are real and simple for $|x| < 1.09\dots \times 10^9$ and whence, by Proposition 2.2, $L_1(x) > 0$ for $|x| < 1.09\dots \times 10^9$. 

Open Problem 4.7 need not be construed as an insurmountable barrier for further research. Indeed, in the interest of new investigations, we propose here a variant of the {\it P\'olyaesque} approach: namely, if you cannot solve a problem change it (for example, generalize it). In this spirit, we mention that
 in the study of the distribution of zeros of entire functions $f(x)\in \s $ (of order $<2$) under the action of the operator $e^{-tD^{2}}$, ($D:=d/dx$) there is a simple heuristic principle formulated by P\'olya.  If $t>0$, then under the action of $e^{-tD^{2}}$ the zeros of $f(x)$ tend to be ``attracted" to the real axis, while under the action of $e^{tD^{2}}$ the zeros of $f$
tend to be repelled by the real axis. Guided by this principle,  we
 apply $e^{-tD^{2}}$ to the Riemann $\xi$-function (see
(4.1)). For convenience and to adhere to the
notation employed in the papers cited below, we
set $H(x):=H_0(x):=\xi(x/2)/8$. Let
$$
H_{t}(x) =e^{-tD^2}H(x)= {\int_{0}^{\infty}}  e^{ts^{2}}
\Phi (s) \cos (xs) ds  \qquad\left ( t \in {\Bbb R};\,  x \in {\Bbb
C},\,\, D:=\frac{d}{dx}\right). \tag4.9
$$  
In 1950, de Bruijn \cite{2} established that (i)  $H_{t} (x)$ has
only real zeros for $t \geq 1/2$ (this is a consequence of the fact
that $H \in\s$, with $\tau=1$, and that $\cos(tD)H \in\lp$ for all
$t \geq 1$) and (ii)  if $H_{t}(x)$ has only
real zeros for some real $t$, then $H_{t'} (x)$ also has only real
zeros for any $t' \geq t$. Subsequently, C. M. Newman \cite{42} showed,  in 1976, that there is a real constant $\Lambda$, which satisfies $-\infty
< \Lambda \leq 1/2$, such that $H_{t}$ has only real zeros if and only
if $t \geq \Lambda$. This constant $\Lambda$ is now called the {\it de
Bruijn-Newman constant} in the literature, and the Riemann Hypothesis
is equivalent to the statement that $\Lambda \leq 0$. Recently, A.M. Odlyzko, W. Smith, R. S. Varga and the author
\cite{17} have shown that $-5.895\cdot 10^{-9} < \Lambda$ (see also \cite{23}).

Differentiation under the integral sign in equation (4.9) (which can be
readily justified by Leibniz's rule, see also \cite{9}) shows that
$H_{t}(x)$ satisfies the backward heat equation:
$$
\frac{\partial (H_{t}(x))}{\partial t} = -
\frac{\partial^{2} H_{t} (x)}{\partial x^{2}} .\tag 4.10
$$
This observation is the key ingredient in the proof of the
following proposition.
\proclaim{ Proposition 4.9}{\rm (\cite{18, Proposition 1})}  Suppose that
$H_{t_{0}}$  has
a multiple real zero. Then $t_{0} \leq \Lambda$.  In particular, if
$t > \Lambda$,  then the zeros of $H_{t}$ {\it are real and simple.}
\endproclaim
\noindent
We next consider two open problems involving $H_{\lambda}(x)$  and the ``new" kernels $\Phi_{\lambda}(t):=e^{\lambda t^2}\Phi(t)$ when (i) $\lambda <0$ and when (ii)  $\lambda >0$. 

\noindent
{\bf Open Problem 4.10.} Fix $\lambda   <0$. Using the theory of positive definite kernels (see Section 3)
 show that  for some non-negative integer $n$, the kernel
$$
K_n(t):=K_n(t; \Phi_{\lambda}):=\int_{-\infty}^{\infty}\Phi_{\lambda}(s+t)\Phi_{\lambda}(s-t)s^{2n}\,ds, \quad\text{is \underbar {not} positive definite.}\tag 4.11
$$
\noindent
Secondly, assume that $\lambda >0$.  In this case, the factor $e^{\lambda s^2}$  under the integral sign (cf. (4.9)) is  an example of a function that P\'olya termed an {\it universal factor} (see the beautiful papers  by P\'olya \cite{50} and de Bruijn \cite{2}).  Universal factors preserve the Laguerre - P\'olya class. In 2009, H. Ki, Y.-O. Kim and J. Lee \cite{36} proved that for every  fixed $\lambda >0$ all but a finite number of the zeros of $H_{\lambda}$ are real and simple. Thus, in particular, if $\lambda >0$, then $H_{\lambda}\in \lp^*$ (see Definition 1.1). Now, in 1987, T. Craven, W. Smith and the author proved the P\'olya-Wiman Conjecture \cite{10} (for a more  elegant proof see  H. Ki and Y.- O. Kim \cite{35}); namely,  if $f(x)\in \lp^*$, then there is a positive integer $m_0$ such that $f^{(m)}(x)\in \lp$ for all $m\ge m_0$. Therefore, it follows from the aforementioned results that for each fixed $\lambda >0$, there is a positive integer $m_0=m_0(\lambda)$ such that for  $2m\ge m_0$ (we work with an even integer so that the new kernel is also even)
$$
H_{\lambda}^{(2m)}(x) =\frac{d^{2m}}{dx^{2m}}e^{-\lambda D^2}H(x)= \int_{0}^{\infty} s^{2m} e^{\lambda s^{2}}
\Phi (s) \cos (xs) ds  \in \lp. \tag4.12
$$  
Observe that the  new kernel $s^{2m}\Phi_{\lambda}(s)=s^{2m} e^{\lambda s^{2}}\Phi (s)$, ($s>0$), is not monotone decreasing, it is not logarithmically concave  and it tends to 0 (as $s\to\infty$) a ``little" slower than $\Phi$.

\noindent
{\bf Open Problem 4.11.} With the above notation and assumptions, is the kernel
$$
K_1(t; \Phi_{\lambda}, m):=\int_{-\infty}^{\infty}\Phi_{\lambda}(s+t)\Phi_{\lambda}(s-t)(s^2-t^2)^m s^2\,ds, \quad\text{positive definite?}
$$

\noindent
We conclude this paper with three additional open problems.

\noindent
{\bf Open Problem 4.12.} Characterize the admissible kernels whose Fourier transforms have all their zeros located in the strip $S(1)$.

\noindent
{\bf Open Problem  4.13.} (\cite{5, Conjecture 2.5}) Show that the derivatives of the Jacobi theta function, $\Phi(t)$, are (strictly) log-concave on  $\Bbb R$;
that is, for each $n\in \Bbb N$,
$$
J_n(t):=(\Phi^{(n)}(t))^2-\Phi^{(n-1)}(t)\Phi^{(n+1)}(t) > 0 \quad\text{for}\quad t\in \Bbb R. \tag 4.13
$$
Since $\Phi(t)$ is an even function (cf. Theorem 4.1), $J_n(t)$ is even and whence it suffices to establish (4.13) for $t\ge 0$.

\noindent Consider again the entire function $F$ (cf. (4.4)) related to the Riemann $\xi$-function: $F(x):=\sum^\infty_{k=0}\frac{ k!\gamma_k}{(2k)!}\frac{x^k}{k!}$, where $\gamma_k:=\frac{k!b_k}{(2k)!}$, ($k=0,1,2\dots$).
Let  $T_k:=\gamma_k^2-\gamma_{k-1}\gamma_{k+1}\ge 0$, (k = 1, 2, 3, \dots), and $E_k:=T_k^2-T_{k-1}T_{k+1}$ for $k=2,3,4,\dots$. Then a necessary condition for the Riemann Hypothesis to hold is that the {\it double Tur\'an inequalities} should hold; i.e., $E_k\ge 0$ for $k=2,3,4\dots$. In \cite{14,  Theorem 2.4}, we derived a concavity condition (for an admissible kernel)  which implies the double Tur\'an inequalities (see also \cite{27, 28}).  Thus, an affirmative answer to the following conjecture, will establish yet another necessary condition for the validity of the Riemann Hypothesis.

\noindent
{\bf Open Problem  4.14.} \cite{14, Problem 3.3} (A new concavity condition of $\Phi(t).$) Let
$s(t):=\Phi(\sqrt{t})$ and set $f(t):=s^\prime(t)^2-s(t)s^{\prime\prime}(t)$. By Theorem 4.2 (b), $f(t)> 0$ for $t>0$.
Then we conjecture that 
$$
\frac{d^2}{dt^2}\log f(t) < 0 \quad \text{for}\ \ t>0. 
$$

\bigskip

\enddocument